\documentclass{article}

\usepackage{amsmath}
\usepackage{amsfonts}
\usepackage{amsthm}
\usepackage{amsopn}
\usepackage{verbatim}

\newcommand{\lb}{\left (}
\newcommand{\rb}{\right )}
\newcommand{\B}[1]{\lb #1 \rb}
\newcommand{\lsb}{\left [}
\newcommand{\rsb}{\right ]}
\newcommand{\SB}[1]{\lsb #1 \rsb}

\newcommand{\be}{\begin{eqnarray}}
\newcommand{\ee}{\end{eqnarray}}
\newcommand{\bc}{}

\DeclareMathOperator{\If}{if}
\DeclareMathOperator{\For}{for}
\DeclareMathOperator{\And1}{and}
\DeclareMathOperator{\that}{that}
\DeclareMathOperator{\is}{is}

\newcommand{\ve}{\varepsilon}

\newcommand{\me}[0]{\mathcal{E}}
\newcommand{\hT}[0]{\hat{T}}
\newcommand{\Bp}[1]{B^{'}\lb #1 \rb}
\newcommand{\Tp}[1][]{T_{+ #1}}
\newcommand{\Tm}[1][]{T_{- #1}}
\newcommand{\hTp}[1][]{\hT_{+ #1}}
\newcommand{\hTm}[1][]{\hT_{- #1}}
\newcommand{\hTpm}[1][]{\hT_{\pm #1}}
\newcommand{\qH}[2][]{q^{#1H_{#2}}}
\newcommand{\leftFactor}[2]{#1[#2] \otimes I}
\newcommand{\rightFactor}[2]{I \otimes #1[#2]}
\newcommand{\leftTwistFactor}[2]{#1[#2] \otimes q^{-2H_{#2}}}
\newcommand{\rightTwistFactor}[2]{q^{2H_{#2}} \otimes #1[#2]}

\newcommand{\YAMLSectionize}[1]{}

\def\dg{\Delta (g)}
\def\gog{g \otimes g}
\def\classlim{q \rightarrow 1}

\newcommand{\Cartan}[1]{\lsb q^{2\phi_i H_i} \rsb}

\def\Chi{\mathcal{X}}

\newcommand{\EChiL}[1][]{\Chi^L_{#1}}
\newcommand{\EChiR}[1][]{\Chi^R_{#1}}
\newcommand{\EChiLT}[1][]{\Chi^{Lt}_{#1}}
\newcommand{\EChiRT}[1][]{\Chi^{Rt}_{#1}}
\newcommand{\EChiD}[1][]{\Chi^\Delta_{#1}}
\newcommand{\EPsiL}[1][]{\Psi^L_{#1}}
\newcommand{\EPsiR}[1][]{\Psi^R_{#1}}
\newcommand{\EPsiLT}[1][]{\Psi^{Lt}_{#1}}
\newcommand{\EPsiRT}[1][]{\Psi^{Rt}_{#1}}
\newcommand{\EPsiD}[1][]{\Psi^\Delta_{#1}}
\newcommand{\QPhiL}[1][]{\Phi^L_{#1}}
\newcommand{\QPhiR}[1][]{\Phi^R_{#1}}
\newcommand{\QPhiD}[1][]{\Phi^\Delta_{#1}}

\newcommand{\matd}[4]{\lb \begin{array}{cc}
#1 & #2 \\ #3 & #4
\end{array} \rb}

\newcommand{\comul}[1]{\Delta \lb #1 \rb}

\newcommand{\delabel}[1]{(\ref{#1})}

\newcommand{\Honed}[1]{\lb \begin{array}{ccc}
#1^{1/2} & 0 \\ 0 & #1^{-1/2}
\end{array} \rb}

\newcommand{\Eoned}[1]{\lb \begin{array}{cc}
1 & #1 \\ 0 & 1
\end{array} \rb}

\newcommand{\Foned}[1]{\lb \begin{array}{cc}
1 & 0 \\ #1 & 1
\end{array} \rb}

\hoffset= - 1.0in         

\title{{\bf The cluster variety face of quantum groups} \vspace{.2cm}}
\author{{\bf A.Popolitov}\thanks{{\small
{\it ITEP, Moscow, Russia and KdVI, University of Amsterdam, the Netherlands}}; popolit@itep.ru, A.Popolitov@uva.nl}}

\begin{document}
 \maketitle

\vspace{-5.0cm}

\begin{center}
\hfill ITEP/TH-05/14\\
\end{center}

\vspace{3.5cm}

\centerline{ABSTRACT}

\bigskip

{\footnotesize
Using the well-known free-field formalism for quantum groups \cite{MV1}, we demonstrate
in case of $A(n)_q$, that quantum group is naturally also a cluster variety \cite{FG1}.
Widely used \cite{Hik1} formulae for mutations \cite{FG1},\cite{FG2} are
direct consequence of independence of group element on the order of simple roots.
Usual formulae \cite{Mars1} for $2 n$ Poisson leaf emerge in classical limit,
if all but few ($2n$) coordinates vanish.
}

\tableofcontents

\bigskip

\bigskip

\section{Introduction}

In modern mathematical physics, one of crucial directions of development is
the search for models, that are integrable, i.e. can be solved exactly.
Having at our disposal sufficiently large toolbox of such models, we
may hope to constructively address real world problems, such as
non-perturbative phenomena in gauge theories.

Traditionally, integrability of a system is
thought to be related to existence of some hidden symmetries, which form a Lie group
(see e.g. \cite{M1}, \cite{M2}, \cite{AASZ} and \cite{LOZ1}).
 If integrable system under consideration is quantum, then
corresponding underlying symmetries also form quantum Lie group.

Recently, a novel way to construct integrable systems, classical as well as quantum,
based on some combinatorial data was presented.
In this approach, discrete integrable systems are governed by so-called cluster algebras \cite{GShTV}, while
continuous integrable systems are governed by cluster varieties \cite{FM}.

Thus, it is natural to ask, whether two kinds of structures, that can control integrability,
are related, and if so, how. In \cite{FG2} an explicit construction of map from cluster variety
to a Lie group was developed, providing a link in one direction. On the other hand,
in \cite{GShV} natural cluster structures on classical Lie groups were studied.

In this paper we show in case of $SL_q(N)$ quantum group, that quantum group possesses natural
structure of quantum $\Chi$-variety of \cite{FG1}. Namely, we restrict ourselves to preferred
family of coordinate systems on quantum group, which is motivated by Gauss's decomposition of
element of classical Lie group. Then equation \eqref{eq:q-group-def},
 defining group subvariety in universal enveloping algebra dictates commutational
relations between coordinate functions, which turn out to have form (49) from \cite{FG1}.
Furthermore, independence of group element from the choice of particular coordinate system
 in the family implies that change from one coordinate system to another
is a substitution of the form (56) from \cite{FG1} (the quantum mutation).

All this poses very intriguing question: can other, less trivial quantum $\Chi$-varieties be
obtained from quantum groups, perhaps, in non-simplylaced cases, or through some reduction procedures?
This, however, is out of scope of the present paper.

This paper is organized as follows: in section \ref{sec:cheatsheet} all the main formulae are
presented without any derivations, mainly for ease of future reference, then
in section \ref{sec:sl2q-case} everything is done for the simplest nontrivial case of $SL_q(2)$,
then in section \ref{sec:slnq-case} general construction for $SL_q(N)$ is presented, and plan
of exposition largely repeats the one in section \ref{sec:sl2q-case}. Finally,
appendix \ref{sec:alpha-beta-gamma-from-chi-psi-phi} is devoted to detailed derivation of commutational
relations for one parametrization of group element from another and appendix \ref{sec:mutation-formulae}
is devoted to some routine checks about validity of mutation formulae.

\section{Cheatsheet}
\label{sec:cheatsheet}

This section is a brief summary of all the important formulae in the paper.
For more detailed explanations, see, for example, next section, where everything is
 done in the simplest case of $SL(2)_q$.

First, we take the quantum group element $g$ to be of the form
of the word of elementary {\it building blocks} $B(i)$, each related to some positive-negative pair
of simple roots $\alpha_{\pm \lsb i \rsb}$
\be
g = B \lb i_1 \rb B \lb i_2 \rb \dots B \lb i_{\frac{n(n + 1)}{2}} \rb,
\ee
where each building block $B(i)$ can be expressed in one of the two equivalent forms
\begin{align}
B(i) & = \me_q \lb \psi_i q^{H_{[i]}} \hT_{+[i]} \rb q^{\phi_i H_{[i]}} \me_{1/q} \lb \chi_i \hT_{-[i]} q^{-H_{[i]}} \rb \\ \notag
& \equiv
w_i^{H_{[i]}} \me_q \lb q^{H_{[i]}} \hT_{+[i]} \rb x_i^{H_{[i]}} \me_{1/q} \lb \hT_{-[i]} q^{-H_{[i]}} \rb y_i^{H_{[i]}}
\end{align}
Here $\me_q(x)$ is $q$-exponential, and $H_{[i]}$, $\hT_{+[i]}$ and $\hT_{-[i]}$ are Chevalley generators of the algebra.
Square brackets denote auxiliary map from all roots to only {\it simple} ones.

Equivalence between two building block parametrizations is given by
\be
\psi_i = w_i,\ \ \chi_i = y_i,\ \ q^{\phi_i} = w_i x_i y_i = y_i x_i w_i
\ee

Quantum group defining equation
\begin{equation}
\label{eq:q-group-def}
\dg = \gog
\end{equation}
 then implies the Darboux-like commutational
relations on the quantum group parameters (also called elements of the dual algebra).
\begin{align}
& \psi_i \chi_j = \chi_j \psi_i, \ \ q^{\phi_i} \psi_j = q^{2 \delta_{ij}} \psi_j q^{\phi_i}, \ \
q^{\phi_i} \chi_j = q^{2 \delta_{ij}} \chi_j q^{\phi_i} & \\ \notag &
y_i w_j = w_j y_i, \ \ x_i w_j = q^{2\delta_{ij}} w_j x_i, \ \ x_i y_j = q^{2\delta_{ij}} y_j x_i
\end{align}

Furthermore, instead of any $B(i)$ one can use slightly different building block $\Bp{i}$
(which also can be written in one of two equivalent forms)
\begin{align}
\Bp{i} & = \me_{1/q} \lb \alpha_i q^{H_{[i]}} \hT_{-[i]} \rb q^{\beta_i H_{[i]}} \me_q \lb \gamma_i \hT_{+[i]} q^{-H_{[i]}} \rb  \\ \notag & \equiv
a_i^{H_{[i]}} \me_{1/q} \lb q^{H_{[i]}} \hT_{-[i]} \rb b_i^{H_{[i]}} \me_q \lb  \hT_{+[i]} q^{-H_{[i]}} \rb c_i^{H_{[i]}},
\end{align}
and {\bf quantum} substitution of variables, which relates thus obtained $g$ to the old one
is of the form of a (quantum) mutation.
\be
a_i = w_i(1 + q x_i), \ \ c_i = y_i(1 + q x_i), \ \ b_i = \frac{1}{x_i},
\ee 
assuming mutation occurred in the building block $B(i)$.

We thus conclude that quantum group $SL(n)_q$ is naturally a cluster  variety,
with very simple cluster data \cite{FG2}:
\begin{itemize}
\item The word of roots $D$, from which seed and frozen seed are constructed is
($i$ denotes $i$-th positive simple root and $\overline{i}$ - negative simple root)
\begin{equation}
D = 1 \overline{1} 2 \overline{2} \dots n \overline{n} 1 \overline{1} 2 \overline{2} \dots (n - 1)\overline{(n-1)}
\dots \dots 1 \overline{1} 2 \overline{2} 1 \overline{1}
\end{equation}
\item $d_i \equiv 1$ for all $i$;
\item $\ve_{ij}$ is almost always zero, only in $x$-$w$- and $x$-$y$-subspaces
  it takes the form of a standard Darboux symplectic matrix, e.g
\begin{equation}
  \ve_{ij} \Big{|}_{x-y-subspace} = \lb \begin{array}{cc} 0 & E \\ -E & 0 \end{array} \rb
\end{equation}

\end{itemize}

\section {$SL(2)_q$ case}
\label{sec:sl2q-case}
In this section we present all parts of the construction in the simplest case - $SL(2)_q$.
Some calculations are shown only in fundamental representation, although
we performed computer checks for higher finite-dimensional representations as well.

\subsection{Quantum algebra's commutational relations and comultiplication rules}
In case of $SL(2)_q$ there are only one positive root $T_{+}$, one negative root $T_{-}$
and one Cartan element $H$.

Commutational relations read
\be
\label{comm-relations-sl2}
q^H \hT_{\pm} = q^{\pm 1} \hT_{\pm} q^H,\ \ \hT_+ \hT_- - \hT_- \hT_+ = \frac{q^{2H} - q^{-2H}}{q - q^{-1}},
\ee
and comultiplication rules are
\be
\Delta(H) = I \otimes H + H \otimes I, \ \ \Delta(\hT_\pm) = q^H \otimes \hT_\pm + \hT_\pm \otimes q^{-H}
\ee

There are two obvious ways to twist generators, associated to roots, which play crucial role in
the construction below.

One is to consider (we call it {\it positive twist})
\be
\label{twist-positively}
T_+ = q^H \hT_+,\ \And1 \ T_- = \hT_- q^{-H}
\ee
with comultiplication rules
\be
\label{comultiply-positively}
\Delta(T_+) = q^{2H} \otimes T_+ + T_+ \otimes I,\ \ \Delta(T_-) = I \otimes T_- + T_- \otimes q^{-2H},
\ee
and another is (we call it {\it negative twist})
\be
\label{twist-negatively}
T_- = q^H \hT_-,\ \And1 \ T_+ = \hT_+ q^{-H},
\ee
with comultiplication rules
\be
\label{comultiply-negatively}
\Delta(T_-) = q^{2H} \otimes T_- + T_- \otimes I,\ \ \Delta(T_+) = I \otimes T_+ + T_+ \otimes q^{-2H},
\ee

Note, that twists preserve commutational relations. That means, that if we are working with only positive-twisted
generators, or with only negative-twisted generators,
we may still choose their fundamental representation to be (like in untwisted case)
\be
T_+ = \matd{0}{1}{0}{0},\ T_- = \matd{0}{0}{1}{0},\ H = \matd{\frac{1}{2}}{0}{0}{-\frac{1}{2}}
\ee

\subsection{Ansatz for group element in two equivalent forms}
In this simple case we have only one pair of positive-negative simple roots, hence group
element contains only one building block
\be
\label{quantum-group-element-sl2}
g = B(1),
\ee
where $B$ can be chosen to be in one of two equivalent forms.

The former is inspired by Fock-Goncharov construction \cite{FG2} and reads
\be
\label{building-block-sl2-fg}
B(1) = & w^H \me_q \lb T_+\rb x^H \me_{1/q} \lb T_- \rb y^H & = \\
= & \Honed{w} \Eoned{1} \Honed{x} \Foned{1} \Honed{y} & = \nonumber \\
= & \matd{w^{1/2}x^{1/2}y^{1/2} + w^{1/2}x^{-1/2}y^{1/2}}{w^{1/2}x^{-1/2}y^{-1/2}}
{w^{-1/2}x^{-1/2}y^{1/2}}{w^{-1/2}x^{-1/2}y^{-1/2}} \nonumber,
\ee

The latter comes from Morozov-Vinet considerations \cite{MV1} and reads
\be
\label{building-block-sl2-mv}
B(1) = & \me_q \lb \psi T_+ \rb q^{\phi H} \me_{1/q} \lb \chi T_- \rb & = \\
= & \Eoned{\psi} \matd{q^{\phi/2}}{0}{0}{q^{-\phi/2}} \Foned{\chi} = \matd{q^{\phi/2} + \psi q^{-\phi/2}\chi}{\psi q^{-\phi/2}}{q^{-\phi/2}\chi}{q^{-\phi/2}} \nonumber.
\ee

In the above formulae $\me_q(x)$ is $q$-exponential (known also as quantum dilogarithm)
\be
\me_q(x) = \sum_{i = 0}^\infty \frac{x^n}{[n]!} q^{-\frac{1}{2}n(n-1)},
\ee
and $q$-analogs are defined to be symmetrical
$$[n] = \frac{q^n - q^{-n}}{q - q^{-1}}$$
Note also, that we have explicitly used the fact, that in fundamental representation $T_+^2 = T_-^2 = 0$,
and $q$-exponential is hence quite simple.

For the reasons, that will become clear in section \ref{comm-relations-dual-algebra-sl2},
$T_+$ and $T_-$ in these formulae should be twisted like in \delabel{twist-positively}
(otherwise there won't be any solution to quantum group defining equation \delabel{eq:q-group-def}).
The other, negative, choice of twist is needed, when we consider building block $\Bp{1}$.

If variables $w$, $x$, $y$ and $\psi$, $\chi$, $q^\phi$ were commutative, then we could immediately
conclude (by equating two different forms of $B(1)$ \delabel{building-block-sl2-fg} and \delabel{building-block-sl2-mv}),
 that correct substitution of variables is
\begin{equation}
  \boxed {
    \label{mv-fg-relation-sl2}
    q^\phi = w x y, \ \ \psi = w,\ \ \chi = y,
  }
\end{equation}

but in quantum (non-commutative) case, everything is not so straightforward.
We have to assume specific form of commutational relations.

Namely, let's assume, that
\begin{equation}
  \label{comm-assumption-sl2}
  \boxed {
    \psi \chi = \chi \psi,\ q^\phi \psi = q^n \psi q^\phi,\ q^\phi \chi = q^n \chi q^\phi,\ \And1 \
    w y = y w,\ x w = q^n w x, \ x y = q^n y x
  }
\end{equation}

for some parameter $n$.
Then
\be
w x y = y x w,\ \And1 \ \lb w^{1/2} x^{1/2} y^{1/2} \rb \lb w^{1/2} x^{1/2} y^{1/2} \rb = w x y,
\ee
and formula \delabel{mv-fg-relation-sl2} is a valid non-commutative substitution of variables.

\subsection{Commutational relations on the dual algebra from comultiplication}
\label{comm-relations-dual-algebra-sl2}

In this section we use quantum group defining equation \delabel{eq:q-group-def} to determine
commutational relations on the dual algebra, following the lines of \cite{MV1}.
 Namely, we take $g$ to be of the form
\delabel{quantum-group-element-sl2}, where $B(1)$ is taken to be of Morozov-Vinet form
\delabel{building-block-sl2-mv}, substitute it into defining equation and provide
such commutational relations for $\phi$, $\psi$ and $\chi$, that the equation becomes valid.




\subsubsection{q-Exponent factorization and comultiplication}
q-Exponent has an interesting property
\be
\me_q(y)\me_q(x) = \me_q(x + y)\ \If \ xy = q^2 yx, \label{fact_q}
\ee
which, for convenience, we will also write in terms of $1/q$ instead of $q$
\be
\me_{1/q}(y)\me_{1/q}(x) = \me_{1/q}(x + y)\ if\ xy = q^{-2} yx, \label{fact_over_q}
\ee

Now comultiplication of the part of the group element corresponding to the positive root factorizes as
(using \delabel{fact_q}; recall, that we use twisting convention \delabel{twist-positively}, hence
comultiplication rules \delabel{comultiply-positively} apply)
\be
\label{comul-expt-positive-sl2}
\comul{\me_q \lb \psi T_+ \rb} & = \me_q \lb \psi \comul{T_+}\rb = \me_q \lb \psi \lb q^{2H} \otimes T_+ + T_+ \otimes I \rb \rb & \\
& = \me_q \lb \psi T_+ \otimes I\rb \me_q \lb \psi q^{2H} \otimes T_+\rb \nonumber
\ee

Analogously, using (\ref{fact_over_q}) we get for the negative root
\be
\label{comul-expt-negative-sl2}
\comul{\me_{1/q} \lb \chi T_- \rb} & = \me_{1/q} \lb \chi \comul{T_-}\rb = \me_{1/q} \lb \chi \lb I \otimes T_- + T_- \otimes q^{-2H} \rb \rb & \\
& = \me_{1/q} \lb \chi T_- \otimes q^{-2H}\rb \me_{1/q} \lb \chi I \otimes T_-\rb \nonumber
\ee

\subsubsection{Convenient notation}

It is convenient to abbreviate frequently occurring combinations of $\me_q$, $T_\pm$, $H$ and $\otimes$.

Namely, let's define the following notation

\be
& \Psi = \me_q \B{\psi \Tp},\ \Chi = \me_{1/q} \B{\chi \Tm},\ \Phi = q^H \nonumber \\
& \EPsiD = \comul{\Psi},\ \EChiD = \comul{\Chi},\ \QPhiD = \comul{\Phi} \nonumber \\
& \EPsiL = \me_q \lb \psi T_+ \otimes I\rb,\ \EPsiR = \me_q \lb \psi I \otimes T_+ \rb \nonumber \\
\label{macrolanguage-sl2}
& \EPsiLT = \me_q \lb \psi T_+ \otimes q^{-2H}\rb,\ \EPsiRT = \me_q \lb \psi q^{2H} \otimes T_+ \rb \\
& \EChiL = \me_{1/q} \lb \chi T_- \otimes I\rb,\ \EChiR = \me_{1/q} \lb \chi I \otimes T_- \rb \nonumber \\
& \EChiLT = \me_{1/q} \lb \chi T_- \otimes q^{-2H}\rb,\ \EChiRT = \me_{1/q} \lb \chi q^{2H} \otimes T_- \rb \nonumber \\
& \QPhiL = q^H \otimes I,\ \QPhiR = I \otimes q^H \nonumber 
\ee

Mnemonic behind this is:
\begin{itemize}
\item $\Psi$, $\Chi$ or $\Phi$ means, that the expression it denotes depends, respectively,
on $\psi$, $\chi$ or $\phi$.
\item $L$ in the superscript means that in the expression there is a tensor product,
and something non-trivial (Chevalley generator) is in the {\it left} part of it.
\item Similarly, $R$ means, that something non-trivial is on the {\it right} of the tensor product entering the expression.
\item Finally, $t$ in the superscript means that tensor product involved is ``twisted'', that is the trivial piece
of the tensor product is multiplied by $q^H$ to a certain power.
\end{itemize}

\subsubsection{Solution of defining equation}

Using the notation \delabel{macrolanguage-sl2}, comultiplication rules \delabel{comul-expt-positive-sl2} and
\delabel{comul-expt-negative-sl2}, together with obvious comultiplication rule for the exponent of Cartan element,
can be expressed as
\be
\EPsiD = \EPsiL \EPsiRT,\ \EChiD = \EChiLT \EChiR,\ \QPhiD = \QPhiL \QPhiR
\ee

Now, left and right hand side of defining equation \delabel{eq:q-group-def} are equal, respectively
\be
& \dg = \EPsiD \QPhiD \EChiD = \EPsiL \EPsiRT \QPhiL \QPhiR \EChiLT \EChiR \label{dg_sl_two} & \\
& \gog = \EPsiL \QPhiL \EChiL \EPsiR \QPhiR \EPsiR \label{gog_sl_two} &
\ee

We immediately see that if
\be
& \EPsiRT \QPhiL = \QPhiL \EPsiR,\ \QPhiR \EChiLT = \EChiL \QPhiR \label{sl_two_first_pass} \\
& \EPsiR \EChiL = \EChiL \EPsiR \label{sl_two_second_pass}
,
\ee
then (\ref{dg_sl_two}) coincides with (\ref{gog_sl_two}).

Equation (\ref{sl_two_second_pass}) is satisfied if
\be
\psi \chi = \chi \psi
\ee

First of (\ref{sl_two_first_pass}) is satisfied if
\be
\psi \lb q^{2H} \otimes T_+ \rb \lb q^{\phi H} \otimes I \rb & = \psi \lb q^{\phi H} \otimes I \rb \lb q^{2H} \otimes T_+ \rb = & \nonumber \\
& = \lb q^{\phi H} \otimes I \rb \psi \lb q^{-2H} \otimes I \rb \lb q^{2H} \otimes T_+ \rb = & \lb q^{\phi H} \otimes I \rb \psi \lb I \otimes T_+ \rb,
\ee

which implies
\be
\psi q^{\phi/2} = q^{\phi/2} \psi q^{-1},\ \that\ \is\ q^\phi \psi = q^2 \psi q^\phi
\ee

Similarly, second of (\ref{sl_two_first_pass}) is satisfied if
\be
q^{\phi H} \chi = \chi q^{\phi H} q^{2H} \\ q^{\phi/2} \chi = \chi q^{\phi/2} q \nonumber \\ q^\phi \chi = q^2 \chi q^\phi \nonumber
\ee

Let's write commutational relations just obtained all on one line

\begin{equation}
\label{comm-relations-psichiphi-sl2}
\boxed{
q^\phi \psi = q^2 \psi q^\phi,\ \ q^\phi \chi = q^2 \chi q^\phi,\ \ \psi \chi = \chi \psi
}
\end{equation}

Thus, our {\it ad hoc} suggestion \delabel{comm-assumption-sl2} proved to be true for commutational
relations on $\psi$, $\chi$ and $\phi$, with $n = 2$.
Now if \delabel{mv-fg-relation-sl2}, understood as {\it quantum} substitution of variables, is correct,
then commutational relations for $x$, $y$ and $w$ are also of the form \delabel{comm-assumption-sl2},
also with $n = 2$
\begin{equation}
\label{comm-relations-wxy-sl2}
\boxed{
x w = q^2 w x,\ \ x y = q^2 y x,\ \ w y = y w
}
\end{equation}
Hence, formula \delabel{mv-fg-relation-sl2} is indeed a full non-commutative change of variables on the
dual algebra of the quantum group $SL_q(2)$.








\subsection{Alternative ansatz for group element}

We can use different form of the building block in the ansatz for the group element.
\be
\label{quantum-group-element-alt-sl2}
g = \Bp{1},
\ee
where building block can be again expressed in two forms
\begin{align}
\label{building-block-alt-sl2-fg}
\Bp{1} & = a^H \me_{1/q} \lb T_-\rb b^H \me_{q} \lb T_+ \rb c^H & \\
& = \Honed{a} \Foned{1} \Honed{b} \Eoned{1} \Honed{c} & \notag \\
& = \matd{a^{1/2}b^{1/2}c^{1/2}}{a^{1/2}b^{1/2}c^{-1/2}}
{a^{-1/2}b^{1/2}c^{1/2}}{a^{-1/2}b^{1/2}c^{-1/2} + a^{-1/2}b^{-1/2}c^{-1/2}} \notag,
\end{align}
or
\begin{align}
\label{building-block-alt-sl2-mv}
\Bp{1} & = \me_q \lb \alpha T_- \rb q^{\beta H} \me_{1/q} \lb \gamma T_+ \rb & \\
& = \Foned{\alpha} \matd{q^{\beta/2}}{0}{0}{q^{-\beta/2}} \Eoned{\gamma}
\\ \notag
& = \matd{q^{\beta/2}}{q^{\beta/2}\gamma}{\alpha q^{\beta/2}}{\alpha q^{\beta/2}\gamma  + q^{-\beta/2}} \notag.
\end{align}

Note, that although the form of $T_\pm$ in this formulae appears to be the same as in
\delabel{building-block-sl2-fg} and \delabel{building-block-sl2-mv}, here they are twisted negatively, like
in \delabel{twist-negatively}. That is, one should be very careful when writing formulae, which include
both types of ansaetze (for example, the ones in Appendix \ref{sec:alpha-beta-gamma-from-chi-psi-phi}, where
we derive commutational relations on $\alpha$, $\beta$ and $\gamma$ from the ones obtained in
\ref{comm-relations-dual-algebra-sl2} for
$\chi$, $\psi$ and $\phi$), as they include additional factors of $q^H$ where appropriate.

Again, naive commutative relation between $(a, b, c)$ and $(\alpha, \beta, \gamma)$ can be
lifted to the full non-commutative relation
\begin{equation}
  \boxed {
    \label{mv-fg-relation-alt-sl2}
    q^\beta = a b c, \ \ \alpha = \frac{1}{a},\ \ \gamma = \frac{1}{c},
  }
\end{equation}
provided the following commutation relations hold
\begin{equation}
  \label{comm-assumption-alt-sl2}
  \boxed {
    \alpha \gamma = \gamma \alpha,\ q^\beta \alpha = q^n \alpha q^\beta,\ q^\beta \gamma = q^n \gamma q^\beta,\ \And1 \
    a c = c a,\ b a = q^{-n} a b, \ b c = q^{-n} c b,
  }
\end{equation}
for some $n$.

\subsubsection{Relation between MV-type ansaetze}
Comparing \delabel{building-block-sl2-mv} to \delabel{building-block-alt-sl2-mv}, one can write
the following four equations
\be
\label{albega-pcp-equations}
& q^{\phi/2} + \psi q^{-\phi/2} \chi = q^{\beta/2} \\
& \psi q^{-\phi/2} = q^{\beta/2} \gamma \nonumber \\
& q^{-\phi/2} \chi = \alpha q^{\beta/2} \nonumber \\
& q^{-\phi/2} = \alpha q^{\beta/2} \gamma + q^{-\beta/2} \nonumber
\ee

First three are sufficient to express $\alpha$, $\beta$ and $\gamma$ through $\chi$, $\psi$ and $\phi$,
while the fourth provides consistency check.

\be
& q^{\beta/2} = q^{\phi/2} + \psi q^{-\phi/2} \chi \nonumber \\
\label{albega-through-pcp}
& \alpha = q^{-\phi/2} \chi q^{-\phi/2} \lb 1 + \psi q^{-\phi/2} \chi q^{-\phi/2} \rb^{-1} \\
& \gamma = \lb 1 + q^{-\phi/2} \psi q^{-\phi/2} \chi \rb^{-1} q^{-\phi/2} \psi q^{-\phi/2} \nonumber
\ee

One can further verify that the form of the commutation relations is preserved by this change of variables, i.e.
\begin{align}
q^{\beta} \alpha = q^2 \alpha q^{\beta} \\ \notag
q^{\beta} \gamma = q^2 \gamma q^{\beta} \\ \notag
\alpha \gamma = \gamma \alpha
\end{align}

The details of the derivation of this commutational relations
can be found in Appendix \ref{sec:alpha-beta-gamma-from-chi-psi-phi}

\subsubsection{Relation between FG-type ansaetze}
Comparing \delabel{building-block-sl2-fg} to \delabel{building-block-alt-sl2-fg} and then
solving the obtained equations, one can express relation between
ansaetze, corresponding to different root orders, in a particularly nice form

\begin{equation}
\label{quantum-mutation-sl2}
  \boxed{
    a = w (1 + q x),\ \ c = y(1 + q x), \ \ b = \frac{1}{x}
  }
\end{equation}
in which one readily recognizes the quantum cluster transformation (also called the mutation).

Taking the limit $\classlim$, of course, reproduces the classical mutation formula
\begin{align}
\label{classical-mutation-sl2}
\notag a = w \lb 1 + x\rb \\
c = y \lb 1 + x\rb \\
b = \frac{1}{x} \notag
\end{align}

The details of the derivation of \delabel{quantum-mutation-sl2} are presented in Appendix
\ref{derivation-of-quantum-mutation-sl2}.




\section {$SL_q(N)$ case}
\label{sec:slnq-case}
In this section describes the construction for arbitrary $A^q_{n} = SL_q(N)$ quantum group.
If you find the discussion too general, you may wish to read $SL_q(2)$ section, and corresponding appendices first.

\subsection{Quantum algebra's commutational relations and comultiplication rules}

Chevalley basis consists of $n$ positive roots, $n$ negative roots and $n$ cartans.
Their (quantum) commutational relations are given by
\be
\label{comm-relations-sln}
\qH{i} \hTpm[j] = q^{\pm C_{ij}/2} \hTpm[j] \qH{i},\ \ 
\hTp[i] \hTm[j] - \hTm[j] \hTp[i] = \frac{\qH[2]{i} - \qH[-2]{i}}{q - q^{-1}},
\ee

There are also commutational relations associated with Serre relations,
but they are left out of the scope of this paper.

Comultiplication rules are as follows
\be
\Delta(H_i) = I \otimes H_i + H_i \otimes I, \ \ 
\Delta(\hT_{\pm i}) = q^{H_i} \otimes \hT_{\pm i} + \hT_{\pm i} \otimes q^{-H_i}
\ee

Again, we may twist each pair of generators $T_{+ i}$ $T_{- i}$ independently in two of the following ways: \\
'positively'
\be
\label{twist-positively-sln}
\Tp[i] = \qH{i} \hTp[i],\ \And1 \ \Tm[i] = \hTm[i] \qH[-]{i},
\ee
or 'negatively'
\be
\label{twist-negatively-sln}
\Tm[i] = \qH{i} \hTm[i],\ \And1 \ \Tp[i] = \hTp[i] \qH[-]{i},
\ee

and the comultiplication rules, correspondingly, read
\be
\label{comultiply-positively-sln}
\Delta(\Tp[i]) = \rightTwistFactor{\Tp}{i} + \leftFactor{\Tp}{i},\ \ 
\Delta(\Tm[i]) = \rightFactor{\Tm}{i} + \leftTwistFactor{\Tm}{i},
\ee
or
\be
\label{comultiply-negatively-sln}
\Delta(\Tp[i]) = \leftTwistFactor{\Tp}{i} + \rightFactor{\Tp}{i},\ \ 
\Delta(\Tm[i]) = \leftFactor{\Tm}{i} + \rightTwistFactor{\Tm}{i},
\ee

If we are working in fundamental representation, and if we restrict ourselves to considering
only 'positively' or only 'negatively' twisted generators, then we may disregard $q^{H_i}$
factors and represent our algebra as follows:
\begin{itemize}
\item $\Tp[i]$ is $N \times N$ matrix of zeroes, except at $i$-th row $i + 1$-th column there is $1$.
\item $\Tm[i]$ is $N \times N$ matrix of zeroes, except at $i+1$-th row $i$-th column there is $1$.
\item $H_i$ is $N \times N$ matrix of zeroes, except at $i$-th row, $i$-th column there is $\frac{1}{2}$
and at $i+1$-th row $i+1$-th column there is $-\frac{1}{2}$
\end{itemize}

\subsection{Ansatz for group element}

In case of $n > 1$ the dimension of the group is larger than the total number of positive, negative
simple roots and cartans.

Hence, in order to express arbitrary group element, we must either use all the roots,
or use simple roots multiple times.

We propose the following, somewhat redundant, parametrization of the group element
\be
\label{quantum-group-element-sln}
g = \prod_{i = 1}^{\frac{n(n+1)}{2}}B(i),
\ee

where each building block $B(i)$ can be in one of two equivalent forms \\
Morozov-Vinet one
\be
\label{building-block-mv-sln}
B(i) = \me_q \B{\psi_i \Tp[\SB{i}]} q^{\phi_i H_{[i]}} \me_{1/q} \B{\chi_i \Tm[\SB{i}]}
\ee

or Fock-Goncharov one
\be
\label{building-block-sln-fg}
B(i) = w_i^{H_{[i]}} \me_q \B{\Tp[\SB{i}]} x_i^{H_{[i]}} \me_{1/q} \B{\Tm[\SB{i}]} y_i^{H_{[i]}},
\ee

where $T_\pm$'s are twisted 'positively' (see \delabel{twist-positively-sln}).

\paragraph{Square bracket map.} Here square brackets around the index $i$ denote, which simple root
we should take for which index. (Clearly, number of indices $\frac{n(n + 1)}{2}$ is larger, than the number
or simple roots, which is $n$)

Namely:
\begin{itemize}
  \item for first $n$ indices we must literally take $i$-th simple root for $i$-th index;
\item for next $n - 1$ indices we must take $(i - n)$-th simple root for $i$-th index;
  \item for yet next $n - 2$ indices we must take $(i - n - (n - 1))$-th simple root for $i$-th index;
\item and so on.
\end{itemize}

Here are the examples of this 'square-bracket' map for first few $n$:
\begin{itemize}
\item n = 1: $\SB{1} = 1$,
\item n = 2: $\SB{1} = 1,\ \SB{2} = 2;\ \
  \SB{3} = 1$,
\item n = 3: $\SB{1} = 1,\ \SB{2} = 2, \SB{3} = 3;\ \ 
  \SB{4} = 1, \SB{5} = 2;
  \SB{6} = 1$
\item n = 4: $\SB{1} = 1,\ \SB{2} = 2, \SB{3} = 3, \SB{4} = 4;\ \
  \SB{5} = 1,\ \SB{6} = 2, \SB{7} = 3;\ \
  \SB{8} = 1,\ \SB{9} = 2;\ \
  \SB{10} = 1$
\end{itemize}


\subsubsection{Relation between MV and FG ansaetze}
If one assumes that $\psi_i$, $\chi_i$ and $\phi_i$ commute in a particular way,
and also $w_i$, $x_i$ and $y_i$ commute in a particular way,
which will turn out to be the case,
then it is easy to write a formula, that connects parameters of MV-ansatz for the building block 
to the parameters of FG-ansatz.

Namely, let's assume, that (note the similarity with \delabel{comm-assumption-sl2} of $SL_q(2)$ case)
\begin{equation}
  \label{comm-assumption-sln}
  \boxed {
    \psi_i \chi_i = \chi_i \psi_i,\ q^{\phi_i} \psi_i = q^n \psi_i q^{\phi_i},\ 
    q^{\phi_i} \chi_i = q^n \chi_i q^\phi_i,\ \And1 \
    w_i y_i = y_i w_i,\ x_i w_i = q^n w_i x_i, \ x_i y_i = q^n y_i x_i
  }
\end{equation}

Then from \delabel{comm-relations-sln}
\be
w_i^{H_{[i]}} \Tp[\SB{i}] = w_i \Tp[\SB{i}] w_i^{H_{[i]}},
\ee
which implies
\be
w_i^{H_{[i]}} \me_q \B{\Tp[\SB{i}]} = \me_q \B{w_i \Tp[\SB{i}]} w_i^{H_{[i]}}
\ee

and similar for $y_i$ and $\Tm[i]$, one gets that
\begin{align}
\label{building-block-fg-transformed-sln}
B(i) & = \me_q \B{w_i \Tp[\SB{i}]} w_i^{H_{[i]}} x_i^{H_{[i]}} y_i^{H_{[i]}} \me_{1/q} \B{y_i \Tm[\SB{i}]}
\\ \notag
& = \me_q \B{w_i \Tp[\SB{i}]} \B{w_i x_i y_i}^{H_{[i]}} \me_{1/q} \B{y_i \Tm[\SB{i}]},
\end{align}
where in the last transition we used the proposed form of commutational relations \delabel{comm-assumption-sln}.

Comparing \delabel{building-block-fg-transformed-sln}
with \delabel{building-block-mv-sln} one clearly sees, that
\begin{equation}
  \boxed {
    \label{mv-fg-relation-sln}
    q^{\phi_i} = w_i x_i y_i, \ \ \psi_i = w_i,\ \ \chi_i = y_i,
  }
\end{equation}
which is a straightforward generalization of \delabel{mv-fg-relation-sl2}. Note, that here we've shown it
to hold in arbitrary representation, not just fundamental.

\subsection{Commutational relations from comultiplication}
Now we must see, that our assumption on the form of commutation relations
\delabel{comm-relations-sln} is indeed true.
So we almost literally reproduce the logic of section \ref{comm-relations-dual-algebra-sl2} here,
only in $SL_q(N)$ case.

\subsubsection{q-Exponent factorization of comultiplication}
Factorization property parallels the one of $SL_q(2)$ case
(see \delabel{comul-expt-positive-sl2}, \delabel{comul-expt-negative-sl2}),
only with indices inserted accordingly
\begin{align}
\label{comul-expt-positive-sln}
\comul{\me_q \B{\psi_i \Tp[\SB{i}]}} =
\me_q \B{\psi_i \leftFactor{\Tp}{\SB{i}}}
\me_q \B{\psi_i \rightTwistFactor{\Tp}{\SB{i}}} \\
\label{comul-expt-negative-sln}
\comul{\me_{1/q} \B{\chi_i \Tm[\SB{i}]}} =
\me_{1/q} \B{\chi_i \leftTwistFactor{\Tm}{\SB{i}}}
\me_{1/q} \B{\chi_i \rightFactor{\Tm}{\SB{i}}}
\end{align}

\subsubsection{Convenient notation}
Again, it is convenient to define a symbolic notation similar to $SL_q(2)$ case, except
everything also has indices now
\be
& \Psi_i = \me_q \B{\psi_i \Tp[\SB{i}]},\ \Chi_i = \me_{1/q} \B{\chi_i \Tm[\SB{i}]},\ \Phi_i = q^{\phi_i H_{\SB{i}}} \nonumber \\
& \EPsiD[i] = \comul{\Psi_i},\ \EChiD = \comul{\Chi_i},\ \QPhiD = \comul{\Phi_i} \nonumber \\
& \EPsiL[i] = \me_q \B{\psi_i \leftFactor{\Tp}{\SB{i}}},\ \EPsiR[i] = \me_q \B{\psi_i \rightFactor{\Tp}{\SB{i}}} \nonumber \\
\label{macrolanguage-sln}
& \EPsiLT[i] = \me_q \B{\psi_i \leftTwistFactor{\Tp}{\SB{i}}},\ 
\EPsiRT[i] = \me_q \B{\psi_i \rightTwistFactor{\Tp}{\SB{i}}} \\
& \EChiL[i] = \me_{1/q} \B{\chi_i \leftFactor{\Tm}{\SB{i}}},\ 
\EChiR[i] = \me_{1/q} \B{\chi_i \rightFactor{\Tm}{\SB{i}}} \nonumber \\
& \EChiLT[i] = \me_{1/q} \B{\chi_i \leftTwistFactor{\Tm}{\SB{i}}},\ 
\EChiRT[i] = \me_{1/q} \B{\chi_i \rightTwistFactor{\Tm}{\SB{i}}} \nonumber \\
& \QPhiL[i] = q^{\phi_i H_{\SB{i}}} \otimes I,\ \QPhiR[i] = I \otimes q^{\phi_i H_{\SB{i}}} \nonumber 
\ee

\subsubsection{Solution of defining equation}
Comultiplication relations \delabel{comul-expt-positive-sln} and \delabel{comul-expt-negative-sln}
 now may be written as
\be
\EPsiD[i] = \EPsiL[i] \EPsiRT[i],\ \EChiD[i] = \EChiLT[i] \EChiR[i],\ \QPhiD[i] = \QPhiL[i] \QPhiR[i]
\ee

Left and right hand side of defining equation \delabel{eq:q-group-def} are equal, respectively
\begin{align}
\dg & = \prod_i \EPsiD[i] \QPhiD[i] \EChiD[i] =
\prod_i \EPsiL[i] \EPsiRT[i] \QPhiL[i] \QPhiR[i] \EChiLT[i] \EChiR[i] \label{dg-sln} & \\
\gog & = \B{\prod_i \EPsiL[i] \QPhiL[i] \EChiL[i]}\B{\prod_i \EPsiR[i] \QPhiR[i] \EPsiR[i]} \label{gog-sln} &
\end{align}

We see, that if
\be
\label{untwist-equations-sln}
& \EPsiRT[i] \QPhiL[i] = \QPhiL[i] \EPsiR[i],\ \And1 \QPhiR[i] \EPsiLT[i] = \EPsiL[i] \QPhiR[i],
\ee

and, furthermore,
\be
\label{left-right-diag-equation-sln}
\EPsiR[i] \EChiL[i] = \EChiL[i] \EPsiR[i],
\ee

then \delabel{dg-sln} can be brought to the form
\begin{align}
\delabel{dg-sln} = & \prod_i \EPsiL[i] \QPhiL[i] \EChiL[i] \EPsiR[i]  \QPhiR[i]  \EChiR[i],
\end{align}

and this expression can be brought to the form \delabel{gog-sln}, provided following holds
\be
\label{left-right-offdiag-equation-sln}
\SB{\Psi\Chi\Phi^L_i,\ \Psi\Chi\Phi^R_j} = 0, \If i \neq j,
\ee

where $\Psi\Chi\Phi$ is {\bf not} a product of 3 expressions, but a {\bf wildcard},
in this place of which any of $\Psi$, $\Chi$ or $\Phi$ can actually stand.

Solving \delabel{untwist-equations-sln} is analogous to solving \delabel{sl_two_first_pass},
one gets
\be
& q^{\phi_i} \psi_i = q^2 \psi_i q^{\phi_i},\ \ q^{\phi_i} \chi_i = q^2 \chi_i q^{\phi_i}
\ee

Solving \delabel{left-right-diag-equation-sln} (cf. solution of \delabel{sl_two_second_pass}) yields
\be
\psi_i \chi_i = \chi_i \psi_i
\ee

Finally, \delabel{left-right-offdiag-equation-sln} is satisfied if
\be
\SB{\psi_i,\ \chi_j} = 0,\ \ \SB{q^{\phi_i},\ \chi_j} = 0,\ \ \SB{q^{\phi_i},\ \psi_j} = 0,\ \For i \neq j
\ee
It is important to note, that since we always commute 'left' exemplars of $\Psi$, $\Chi$, $\Phi$ with
'right' exemplars of $\Psi$, $\Chi$, $\Phi$, additional terms, coming from Serre relations, do not arise.

Thus, commutational relations just obtained can be written as
\begin{equation}
\label{comm-relations-psichiphi-sln}
\boxed{
  q^{\phi_i} \psi_j = q^{2\delta_{ij}} \psi_j q^{\phi_i},\ \ 
  q^{\phi_i} \chi_j = q^{2\delta_{ij}} \chi_j q^{\phi_i},\ \ 
  \psi_i \chi_j = \chi_j \psi_i,
}
\end{equation}

which is consistent with our initial assumption on the form of commutational relations \delabel{comm-assumption-sln}.

From this commutational relations, using the transformation \delabel{mv-fg-relation-sln}
one can show that commutational relations on $w$'s, $x$'s and $y$'s also have the desired form
\begin{equation}
\boxed{
\label{comm-relations-wxy-sln}
w_i y_j = y_j w_i,\ x_i w_j = q^{2 \delta_{ij}} w_j x_i, \ x_i y_j = q^{2 \delta_{ij}} y_j x_i
}
\end{equation}
 thus proving the validity of \delabel{mv-fg-relation-sln} itself, as a non-commutative change of variables.

\subsection{The cluster variety data}

Looking at the commutation relations \delabel{comm-relations-wxy-sln} one can already
recognize the appearance of quantum $\mathcal{X}$-variety structure of \cite{FG1},
and read off the corresponding cluster data.

Namely, vector of integers $d$ is in fact vector of ones, since our group is simply laced.

Then, comparing of ansatz for the group element \delabel{quantum-group-element-sln}
with map to group, presented in section 3 of \cite{FG2}, we see, that
seed $J$ and frozen seed $J_0$ correspond to the following word $D$ of simple roots
\begin{equation}
D = 1 \overline{1} 2 \overline{2} \dots n \overline{n} 1 \overline{1} 2 \overline{2} \dots (n - 1)\overline{(n-1)}
\dots \dots 1 \overline{1} 2 \overline{2} 1 \overline{1},
\end{equation}
where $i$ denotes $i$-th positive simple root, and $\overline{i}$ - corresponding negative simple root.

Finally, the skew-symmetric matrix $\ve$ is nonzero only in x-y and x-w subspaces and
there takes the form of a standard symplectic form
\begin{equation}
  \ve_{ij} \Big{|}_{x-y-subspace} = \ve_{ij} \Big{|}_{x-w-subspace} = \lb \begin{array}{cc} 0 & E \\ -E & 0 \end{array} \rb
\end{equation}


\subsection{Relation between FG-type ansaetze}

Analogously to $SL_q(2)$ case, one can consider different parametrization
of a building block $\Bp{i}$

\begin{align}
\label{building-block-alt-sln-fg}
\Bp{i} & = a_i^{H_{\SB{i}}} \me_{1/q} \lb \Tm[\SB{i}] \rb b_i^{H_{\SB{i}}} \me_{q} \lb \Tp[\SB{i}] \rb c_i^{H_{\SB{i}}} &
\end{align}

In order for \delabel{building-block-alt-sln-fg} to be equal to the other parametrization
 \delabel{building-block-sln-fg}, we must have the following relation between parameters

\begin{equation}
a_i = w_i (1 + q x_i),\ \ c_i = y_i ( 1 + q x_i), \ \ b_i = \frac{1}{x_i}
\end{equation}

And again, for this to be valid, the following equation should hold

\begin{align}
\label{mutation-sln-main-equation}
\me_q(q^{H_{\SB{i}}} \hT_{+ \SB{i}})
x_i^{H_{\SB{i}}}
\me_{1/q}(\hT_{- \SB{i}} q^{-H_{\SB{i}}})
= \me_q\lb \frac{q^{2H_{\SB{i}}} x_i}{q - 1/q}\rb \me_{1/q}(q^{H_{\SB{i}}} \hT_{-\SB{i}})
x_i^{-H_{\SB{i}}}
\me_{q}(\hT_{+ \SB{i}} q^{-H_{\SB{i}}}) \me_{1/q}\lb \frac{q^{-2H_{\SB{i}}} x_i}{1/q - q}\rb
\end{align}

Since, we change only the variables associated with $i$-th index,
 essentially we are reparametrizing some $SL_q(2)$ subgroup of $SL_q(N)$,
so for details see appendix \ref{sec:mutation-formulae}.

Note, that here we've in fact considered only simplest possible changes
in FG parametrization. In FG construction it is also possible to perform mutations
in other directions. Roughly speaking, these other mutations should correspond to
to the interchange of $\Tp{\SB{i}}$ with $\Tm{\SB{i}}$, that are not in the same building block $B(i)$.
 However, it is now not clear how to precisely do that, so this is the subject of
further investigation.

\subsection{Reduction to $2n$ symplectic leaf}

It is clear, that if we consider only subalgebra of functions of first $n$ $x_i$'s and first $n$ $w_i$'s,
it is closed under commutational relations \delabel{comm-relations-wxy-sln}. This means,
that it defines a $2n$-dimensional submanifold in quantum group $SL_q(N)$.
Points of this manifold can be explicitly parametrized as

\begin{equation}
g_{2n} = \prod_{i = 1}^n w_i^{H_{i}} \me_q \B{\Tp[i]} x_i^{H_{i}} \me_{1/q} \B{\Tm[i]}
\end{equation}

Taking the limit $\classlim$ in this expression, one readily recognizes
parametrization of $2n$ symplectic manifold of a Lie group from \cite{Mars1},
on which a relativistic Toda chain lives.

Thus, it is natural to think, that on $g_{2n}$ a $q$-deformation of relativistic Toda chain
system lives.

\section{Conclusion}
In this paper we looked at $SL_q(2)$ and $SL_q(N)$ quantum groups through
the lens of Morozov-Vinet construction of free-field representation of
quantum group element \cite{MV1}. Namely, by a clever choice of order of factors, we
ensured that classical limit $\classlim$ of a group element $g$ looks like
Fock-Goncharov map of a cluster variety into a Lie group \cite{FG2}.
We demonstrated, that some different choices of MV parametrization are equivalent
and related by quantum cluster mutations.
We also showed, how to make a reduction to a quantum version of $2n$ Poisson submanifold,
on which a natural integrable system thus lives \cite{Mars1}.
We interpret our results as strongly indicating that quantum group $SL_q(N)$ naturally
caries the structure of a quantum cluster variety, in fact with very simple cluster data.

Also, as an experiment, all versions of LaTeX source of this article, as well as
photos of some drafts of calculations, used in preparation
and some Mathematica plain-text files (conveniently usable through corresponding Emacs mode), will be available as
a publicly read-accessible Git repository, at
\\ https://github.com/mabragor/quantum-group-looks.git

\section*{Acknowledgments}

Author is indebted to A.Mironov and A.Morozov for numerous stimulating discussions.
Author also kindly thanks ETH Zurich and GQT national Dutch mathematical research cluster,
as parts of work were done during their conferences and schools.
Our work is partly supported by RFBR grants 14-01-31492\_mol\_a and 14-02-00627,
 grant for support of scientific schools NSh-1500.2014.2 and Vici grant of the NWO.

\appendix

\section{Commutational relations on $(\alpha\ \beta\ \gamma)$ from that of $(\chi\ \psi\ \phi)$}
\label{sec:alpha-beta-gamma-from-chi-psi-phi}

Recall, that $(\alpha\ \beta\ \gamma)$ is expressed through $(\chi\ \psi\ \phi)$ via formula \delabel{albega-through-pcp}
\be
& q^{\beta/2} = q^{\phi/2} + \psi q^{-\phi/2} \chi \nonumber \\
& \alpha = q^{-\phi/2} \chi q^{-\phi/2} \lb 1 + \psi q^{-\phi/2} \chi q^{-\phi/2} \rb^{-1} \nonumber \\
& \gamma = \lb 1 + q^{-\phi/2} \psi q^{-\phi/2} \chi \rb^{-1} q^{-\phi/2} \psi q^{-\phi/2} \nonumber 
\ee

First, we want to check, that, given this, the fourth equation of \delabel{albega-pcp-equations} also holds.
For this we take as input, that commutational relations \delabel{comm-relations-psichiphi-sl2} hold.
\\
Indeed
\def\qbHalf{\lb q^{\phi/2} + \psi q^{-\phi/2} \chi \rb}
\def\qbMinusHalf{\qbHalf^{-1}}
\def\qbHalfNorm{\lb 1 + q^{-\phi/2} \psi q^{-\phi/2} \chi \rb}
\def\qbMinusHalfNorm{\qbHalfNorm^{-1}}
\be
\alpha q^{\beta/2} \gamma + q^{-\beta/2}
= & q^{-\phi/2}\chi \qbMinusHalf \psi q^{-\phi/2}
+ \qbMinusHalf & = \nonumber \\
\label{albega-pcp-derivation-end1}
= & q^{-\phi/2} \chi \qbMinusHalfNorm q^{-\phi/2} \psi q^{-\phi/2} + \qbMinusHalfNorm q^{-\phi/2} & =
\ee

To transform the expression further, we need two auxiliary elementary facts
\be
& \chi \cdot \lb q^{-\phi/2} \psi q^{-\phi/2} \chi \rb
= q^{2} \lb q^{-\phi/2}\psi q^{-\phi/2} \chi \rb \cdot \chi \nonumber \\
& q^{-\phi/2} \cdot \lb q^{-\phi/2} \psi q^{-\phi/2} \chi \rb
= q^{-2} \lb q^{-\phi/2}\psi q^{-\phi/2} \chi \rb \cdot q^{-\phi/2} \nonumber
\ee

So we get (this can be seen by representing $\qbMinusHalfNorm$ as infinite series and commuting $q^{-\phi/2}\chi$
with all individual monomials
\be
\delabel{albega-pcp-derivation-end1} = & \qbMinusHalfNorm q^{-\phi/2} \chi q^{-\phi/2} \psi q^{-\phi/2}
+ \qbMinusHalfNorm q^{-\phi/2} & = \nonumber \\
= & \qbMinusHalfNorm \lsb q^{-\phi/2} \chi q^{-\phi/2} \psi q^{-\phi/2} + q^{-\phi/2}\rsb = \label{albega-pcp-derivation-end2}
\ee

The first summand in the square brackets can be simplified as
\be
q^{-\phi/2} \chi q^{-\phi/2} \psi q^{-\phi/2} = q^{-\phi/2} q q^{-\phi/2} \chi \psi q^{-\phi/2}
= q^{-\phi/2} q q^{-\phi/2} \psi \chi q^{-\phi/2} = q^{-\phi/2} \psi q^{-\phi/2} \chi q^{-\phi/2} \nonumber
\ee

and we finally get
\be
\delabel{albega-pcp-derivation-end2} = & \qbMinusHalfNorm \qbHalfNorm q^{-\phi/2} = q^{-\phi/2},
\ee
which completes the check.

Second, we want to derive commutational relations on $(\alpha\ \beta\ \gamma)$ from commutational
relations on $(\chi\ \psi\ \phi)$.

First, note that if
\be
\lb q^{\beta/2}\gamma \rb \cdot = q^{-n} q^{\beta/2} \cdot \lb q^{\beta/2}\gamma \rb,
\ee
for some $n$, then also
\be
 \gamma \cdot q^{\beta/2} = q^{-n} q^{\beta/2} \cdot \gamma,
\ee
and analogous statement is valid for $\alpha$
 (except we consider $\alpha q^{\beta/2}$ in place of $q^{\beta/2} \gamma$) 

Then, using explicit expressions \delabel{albega-through-pcp} we may write
\be
q^{\beta/2} \gamma \cdot q^{\beta/2} = & \psi q^{-\phi/2} \qbHalf
= q^{-1} q^{\phi/2} \psi q^{-\phi/2} + \psi q^{-2} \psi q^{-\phi/2} \chi q^{-\phi/2} & = \\
& = q^{-1} q^{\phi/2} \psi q^{-\phi/2} + q^{-1} \psi q^{-\phi/2} \chi \psi q^{-\phi/2} = q^{-1} q^{\beta/2} \cdot q^{\beta/2} \gamma,
\ee
so it turned out that $n = 1$. And ditto for $\alpha$
\be
\alpha q^{\beta/2} \cdot q^{\beta/2} = & q^{-\phi/2} \chi \qbHalf
= q^{-1} q^{\phi/2} q^{-\phi/2} \chi + q^{-\phi/2} q \psi q^{-\phi/2} \chi \chi & = \\
& = q^{-1} q^{\phi/2} q^{-\phi/2} \chi + q^{-1} \psi q^{-\phi/2} \chi q^{-\phi/2} \chi
= q^{-1} q^{\beta/2} \cdot \alpha q^{\beta/2},
\ee
so $n$ here also equals $1$.

The only thing that remains to be shown, is that $\alpha$ and $\gamma$ commute.
For that, let's note that if
\be
\alpha q^{\beta/2} \cdot q^{\beta/2} \gamma = q^m q^{\beta/2} \gamma \alpha q^{\beta/2},
\ee
for some $m$, then
\be
\alpha \gamma = q^m \gamma \alpha
\ee
Indeed,
\be
\alpha \gamma = \alpha q^{-\beta/2} q^{\beta/2} q^{\beta/2} q^{-\beta/2} \gamma
= q^n q^{-n} q^{-\beta/2} \alpha q^{\beta/2} q^{\beta/2} \gamma q^{-\beta/2}
= q^m q^{-\beta/2} q^{\beta/2} \gamma \alpha q^{\beta/2} q^{-\beta/2} = q^m \gamma \alpha
\ee

Now, one easily sees, that
\be
\alpha q^{\beta/2} q^{\beta/2} \gamma = q^{-\phi/2} \chi \psi q^{-\phi/2} = q q^{-\phi/2} \psi q^{-\phi/2} \chi
= \psi q^{-\phi/2} q^{-\phi/2} \chi = q^{\beta/2} \gamma \alpha q^{\beta/2},
\ee
so $m = 0$.

\section{Mutation formulae}
\label{sec:mutation-formulae}
\subsection{$SL_q(2)$ fundamental representation}
\label{derivation-of-quantum-mutation-sl2}
This derivation is straightforward, yet it is example of calculation,
which already contains some peculiarities, which typically arise when doing non-commutative algebra,
hence we present it here in very detailed form.

We start from the following equations, which come from equating \delabel{building-block-sl2-fg} to \delabel{building-block-alt-sl2-fg}
\be
& w^{1/2}x^{1/2} y^{1/2} + w^{1/2}x^{-1/2}y^{1/2} = a^{1/2}b^{1/2}c^{1/2} \\
& w^{1/2}x^{-1/2}y^{-1/2} = a^{1/2}b^{1/2}c^{-1/2} \\
& w^{-1/2}x^{-1/2}y^{1/2} = a^{-1/2}b^{1/2}c^{1/2} \\
& w^{-1/2}x^{-1/2}y^{-1/2} = a^{-1/2}b^{1/2}c^{-1/2} + a^{-1/2}b^{-1/2}c^{-1/2}
\ee
From second equation we get (multiplying both sides by $c$ to the right)
\be
& w^{1/2}x^{-1/2}y^{-1/2} c = a^{1/2}b^{1/2}c^{1/2},
\ee
which, using the first equation allows us to state that
\be
c = y^{1/2}x^{1/2}w^{-1/2}\lb w^{1/2}x^{1/2}y^{1/2} + w^{1/2}x^{-1/2}y^{1/2} \rb
= y^{1/2} x y^{1/2} + y = y\lb 1 + qx\rb = \lb 1 + \frac{x}{q} \rb y,
\ee
where we used commutational relations \delabel{comm-relations-wxy-sl2}.

Analogously, multiplying both sides of the third equation by a to the left, we get
\be
& a w^{-1/2} x^{-1/2} y^{1/2} = a^{1/2}b^{1/2}c^{1/2},
\ee
which, using the first equation leads
\be
a = \lb w^{1/2}x^{1/2}y^{1/2} + w^{1/2}x^{-1/2}y^{1/2} \rb y^{-1/2}x^{1/2}w^{1/2}
= w^{1/2}x w^{1/2} + w = w \lb 1 + q x \rb
\ee

Now the trickiest part is to calculate the expression for $b$. Naively, one  would want
to multiply the first equation by $a^{-1/2}$ to the left and by $c^{-1/2}$ to the right,
get $b^{1/2}$ and then square it.
However, calculating $a^{1/2}$ is not as simple as taking the product of square roots of both factors
\be
&\sqrt{a} \neq \sqrt{w} \sqrt{1 + q x},
\ee
instead there is this q-Exponential formula \delabel{a-h-through-w}, which one will need when
doing calculation in arbitrary representation.
Failure to notice this peculiarity easily leads to incorrect expression for $b$
\be
b \neq \frac{\lb 1 + x\rb^2}{\lb 1 + q x\rb\lb 1 + \frac{x}{q}\rb} \frac{1}{x}.
\ee
Here, on the other hand, we can simply observe, that
\be
\sqrt{a} \sqrt{b} \sqrt{c} \sqrt{a} \sqrt{b} \sqrt{c} = q^{-1/2} a \sqrt{b} \sqrt{c} \sqrt{b} \sqrt{c} = abc,
\ee
and then, using the first equation

\be
b = & a^{-1} \lb \sqrt{a}\sqrt{b}\sqrt{c} \rb^2 c^{-1} & \nonumber \\
= & \lb 1 + qx\rb^{-1}w^{-1}\lsb w^{1/2}x^{1/2}y^{1/2} + w^{1/2}x^{-1/2}y^{1/2} \rsb
\lsb w^{1/2}x^{1/2}y^{1/2} + w^{1/2}x^{-1/2}y^{1/2} \rsb y^{-1}\lb 1 + \frac{x}{q}\rb^{-1}
\\ = & \lb 1 + qx\rb^{-1} \lsb x + \lb q + \frac{1}{q} \rb \frac{1}{x} \rsb \lb 1 + \frac{x}{q}\rb^{-1} = \frac{1}{x}
\nonumber 
\ee

The check, that the fourth equation is valid is trivial, we do not present it here.

\subsection{$SL_q(2)$ general case}
\paragraph{Formula for $a^H$}
As was mentioned above, one should be careful, when writing down formula for $a^n$ in terms of $w$.

First few examples
\begin{align}
& a^1 = w(1 + q x) & \\ \notag
& a^2 = w(1 + q x) w (1 + q x) = w^2 (1 + q^3 x) (1 + q x) & \\ \notag
& a^3 = w(1 + q x) w (1 + q x) w (1 + q x) = w^3 (1 + q^5 x) (1 + q^3 x) (1 + q x) &
\end{align}

and in general we get
\begin{align}
& a^n = w^n \prod_{i = 1}^n (1 + q^{2 i - 1} x) = w^n
  \frac{\prod_{i = 1}^\infty (1 + q^{2i -1} x)}{\prod_{i = 1}^\infty (1 + q^{2n} q^{2i -1} x)}
\end{align}

and the last expression can be easily generalized for matrix-valued $n$'s.
We are specifically interested in $n = H$
\begin{align}
\label{a-h-through-w}
a^H & = w^H \frac{\prod_{i = 1}^\infty (1 + q^{2i -1} x)}{\prod_{i = 1}^\infty (1 + q^{2H} q^{2i -1} x)} & \\ \notag
& = w^H \frac{\me_q \lb q^{2H} \frac{x}{q - 1/q} \rb}{\me_q \lb \frac{x}{q - 1/q} \rb}
\end{align}

Analogously, for $c^H$ we get (this time we push all $y$'s to the right)
\begin{align}
\label{c-h-through-y}
c^H = \frac{\me_{1/q} \lb q^{-2H} \frac{x}{q - 1/q} \rb}{\me_{1/q} \lb \frac{x}{q - 1/q} \rb} y^H
\end{align}

Obviously,
\begin{align}
\label{b-h-through-x}
& b^H = x^{-H}
\end{align}

\paragraph{Mutation equation}
Again, equating \delabel{building-block-sl2-fg} to \delabel{building-block-alt-sl2-fg},
and using formulae \delabel{a-h-through-w}, \delabel{c-h-through-y} and \delabel{b-h-through-x}
we arrive at
\begin{align}
\label{mutation-sl2-main-equation}
\me_q(q^H \hT_+) x^H \me_{1/q}(\hT_- q^{-H}) = \me_q\lb \frac{q^{2H} x}{q - 1/q}\rb
\me_{1/q}(q^H \hT_-) x^{-H} \me_{q}(\hT_+ q^{-H}) \me_{1/q}\lb \frac{q^{-2H} x}{1/q - q}\rb
\end{align}

Note, that $T_\pm$'s were expressed through $\hT_\pm$'s differently at the l.h.s and the r.h.s.

We've checked in Mathematica explicitly, that this equation holds for first 20 symmetric representations,
but the general proof is still missing.

\paragraph{Remarks on infinite-dimensional representation}
Here we provide example of check of mutation equation for infinite dimensional representation.

Generators $H$ and $\hT_\pm$'s can be represented as differential operators,
either acting on space of functions of $x$ and $y$;
\begin{equation}
\label{eq:sl2q-diff-hom-repr}
H = \frac{1}{2} x \partial_x - \frac{1}{2} y \partial_y,\ \ \ \hT_+ = x \partial^q_y,\ \ \ \hT_- = y \partial^q_x,
\end{equation}
or as differential operators, acting on functions of $x$ only (this time there is 1-parametric
family of possible representations, parametrized by $j$)
\begin{equation}
\label{eq:sl2q-diff-nonhom-repr}
H = x \partial_x - j,\ \ \ \hT_+ = [2 j] x M_q - q^{2j} x^2 \partial^q_x,\ \ \ \hT_- = \partial^q_x
\end{equation}
Here $\partial^q_x$ is a (symmetric) $q$-derivative w.r.t $x$ (analogously for $y$) and
$M_q$ is operator of dilation of $x$
\begin{equation}
\partial^q_x f(x) = \frac{f(qx) - f(x / q)}{q - 1/q},\ \ \ M_q f(x) = f(q x)
\end{equation}

We will consider representation, corresponding to highest weight $x^{-2}$ in both representations
(in second representation we then must put $j = -1$).

After that, in basis $x^{-2 + k}$, $k = 0, 1, 2, \dots$ generators become semi-infinite matrices
\begin{equation}
H = \lb \begin{array}{cccc}
  -1 & 0 & 0 & . \\
  0 & -2 & 0 & . \\
  0 & 0 & -3 & . \\
  . & . & . & \dots \\
\end{array}\rb,\ \ \ 
\hT_+ = \lb \begin{array}{ccccc}
  0 & [1] & 0 & 0 & . \\
  0 & 0 & [2] & 0 & . \\
  0 & 0 & 0 & [3] & . \\
  0 & 0 & 0 & 0 & . \\
  . & . & . & . & \dots \\
\end{array}\rb,\ \ \ 
\hT_- = \lb \begin{array}{ccccc}
  0 & 0 & 0 & 0 & . \\
  -[2] & 0 & 0 & 0 & . \\
  0 & -[3] & 0 & 0 & . \\
  0 & 0 & -[4] & 0 & . \\
  . & . & . & . & \dots \\
\end{array}\rb,\ \ \ 
\end{equation}
and it is straightforward to calculate q-Exponential building blocks in equation \eqref{mutation-sl2-main-equation}
(we write formulae for element in $i$-th row $j$-th column, when it's nonzero)

\newcommand{\qbin}[2]{\lsb \begin{array}{c} #1 \\ #2 \end{array} \rsb}

\begin{align}
\me_q(q^H \hT_+)_{ij} = \qbin{j - 1}{i - 1} q^{-(j-i)(j-1)},\ \ \
\me_{1/q}(\hT_- q^{-H})_{ij} = \qbin{i}{j} (-1)^{i - j} q^{(i-j)(i-1)} \\ \notag
\me_{1/q}(q^H \hT_-)_{ij} = \qbin{i}{j} (-1)^{i-j} q^{-(i-j)(j+1)},\ \ \
\me_{q}(\hT_+ q^{-H})_{ij} = \qbin{j-1}{i-1} q^{(j-i)(i+1)},
\end{align}
where $\qbin{n}{k}$ denotes quantum binomial coefficient.

Then, $n$-$m$-th matrix element of equation \eqref{mutation-sl2-main-equation}
, in region where $m \geq n$,
is the following equality between hypergeometric functions (for simplicity we provide formula for $q = 1$)
\begin{equation}
\label{eq:mutation-equation-sl2-corollary}
x^{-m} \frac{\Gamma(m)}{\Gamma(n) \Gamma(m - n + 1)} {_2F_1} \lb m, m+1 ; m - n + 1; -\frac{1}{x}\rb
= \frac{(-1)^{n+1} n x}{(1 + x)^{m + n}} {_2F_1} \lb 1-m, 1-n ; 2; -x \rb
\end{equation}

If we put $m = n + k$, for each particular choice of $n$ this becomes elementary function in $x$ and $k$,
so we checked, that the equality \eqref{eq:mutation-equation-sl2-corollary} holds for first 25 $n$ and
arbitrary $x$ and $k$.

Generalization to $q \neq 1$ seems straightforward (just substitute all functions by their $q$-analogs),
however, we didn't perform any checks in this case.

\end{document}